\documentclass[12pt]{article}
\usepackage{authblk}
\usepackage{amsmath}
\usepackage[english]{babel}
\usepackage{theorem}
\usepackage[T1]{fontenc}
\usepackage{amssymb,xspace}
\usepackage{tikz}
\usepackage[applemac]{inputenc}

\newcommand{\Jeu}{{\sc influence}\xspace}

\usepackage{theorem}
\newtheorem{theorem}{Theorem}
\newtheorem{lemma}[theorem]{Lemma}
\newtheorem{corollary}[theorem]{Corollary}
\newtheorem{proposition}[theorem]{Proposition}
\newtheorem{definition}[theorem]{Definition}
\newtheorem{conjecture}[theorem]{Conjecture}
\newtheorem{ex}[theorem]{Example}

\tikzstyle{vertex}=[circle, fill=black, inner sep= 0, minimum size = 6]
\tikzstyle{vertexW}=[circle, draw, fill=white, inner sep= 0, minimum size = 6]
\tikzstyle{edge}=[->, thick]

\newtheorem{theo_bis}[theorem]{\noindent\bf Theorem}%
{\begin{theo_bis}{\bf #1} \sl #2}%
{\end{theo_bis}}

\newtheorem{propo_bis}[theorem]{\noindent\bf Proposition}%
{\begin{propo_bis}{\bf #1} \sl #2}%
{\end{propo_bis}}

\newtheorem{r_bis}[theorem]{\noindent\bf Remark}
\newenvironment{remark}[2][]%
{\begin{r_bis}{\bf #1} \normalfont #2}%
{\hfill $\Box$\end{r_bis}}

\newtheorem{col_bis}[theorem]{\noindent\bf Corollary}%
{\begin{col_bis}{\bf #1} \sl #2}%
{\end{col_bis}}

\newenvironment{preuve}[2][]
{\begin{trivlist}{\bf\item Proof} \textbf{#1} #2}
{\hfill $\blacksquare$\end{trivlist}}

\long\def\symbolfootnote[#1]#2{\begingroup\def\thefootnote{\fnsymbol{footnote}}\footnote[#1]{#2}\endgroup} 

{\vspace{-0.2cm}\begin{list}%
	{--}%
	{\setlength{\labelwidth}{30pt}%
	 \setlength{\leftmargin}{35pt}%
	 \setlength{\itemsep}{-0.1cm}}}%
{ \end{list}\vspace{-0.2cm} }

\title{\Jeu : a partizan scoring game on graphs\footnote{This research has benefited from the financial support of IDEXLYON from Universit\'e de Lyon (project INDEPTH) within the Programme Investissements d'Avenir (ANR-16-IDEX-0005) and the SFR ``Math\'ematiques de la d\'ecision pour l'ing\'enierie physique et sociale'' (MODMAD).}}
\author[1]{Eric Duch\^ene}
\author[2]{Stphane Gonzalez}
\author[1]{Aline Parreau}
\author[2]{Eric R\'emila}
\author[2]{Philippe Solal}

\affil[1]{Universit\'e de Lyon, Universit\'e Lyon 1, LIRIS UMR CNRS 5205\\

F-69621, Lyon, France}
\affil[2]{Universit\'e de Lyon, UJM St-Etienne,  GATE Lyon Saint-Etienne UMR CNRS 5824, F-42023, St-Etienne, France}

\begin{document}

\maketitle

\begin{abstract}
  We introduce the game \Jeu, a scoring combinatorial game, played on a directed graph where each vertex is either colored black or white. The two players, Black and White play alternately by taking a vertex of their color and all its successors (for Black) or all its predecessors (for White). The score of each player is the number of vertices he has taken.
  We prove that \Jeu is a nonzugzwang game, meaning that no player has interest to pass at any step of the game, and thus belongs to Milnor's universe. We study this game in the particular class of paths where black and white are alternated. We give an almost tight strategy for both players when there is one path. More precisely, we prove that the first player always gets a strictly better score than the second one, but that the difference between the score is bounded by $5$. Finally, we exhibit some graphs for which the initial proportion of vertices of the color of a player is as small as possible but where this player can get almost all the vertices.
\end{abstract}

\section{Introduction}


We define the following $2$-player game, namely \Jeu, on a directed graph where each vertex is either colored black or white. The two players, namely Black and White (or simply $B$ and $W$ for short), play alternately by choosing a vertex $v$ of their own color and:

\begin{itemize}
  \item if it is Black's turn, remove $v$ and all the vertices reachable from $v$ through a directed path;
  \item if it is White's turn, remove $v$ and all the vertices for which there exists a  directed path to $v$.
\end{itemize}

The game ends when all the vertices have been removed. When a player has no more vertex of his own color, he has to wait for his opponent to empty the directed graph. The score of each player is the number of vertices he has removed. Each player wants to maximize his score. \\

\begin{figure}[h]
			\centering
			\begin{tikzpicture}
			\node[vertex] (0) at (0,2) {};
			\node[vertexW] (1) at (2,2) {};
			\node[vertex] (2) at (1,1) {};
			\node[vertexW] (3) at (3,1) {};
			\node[vertexW] (4) at (0,0) {};
			\node[vertexW] (5) at (2,0) {};
			\draw[edge](0)--(2);
			\draw[edge](1)--(2);
			\draw[edge](2)--(4);
			\draw[edge](2)--(5);
			\draw[edge](3)--(5);
			\draw[edge](3)--(1);
            \draw(0) node[left=5] {$u$};
             \draw(1) node[right=5] {$v$};
              \draw(2) node[left=5] {$w$};
               \draw(3) node[right=5] {$x$};
                \draw(4) node[left=5] {$y$};
                 \draw(5) node[right=5] {$z$};
			
			\end{tikzpicture}
			
	\caption{A game position where the first player wins}
	\label{fig-ex1}
	\end{figure}
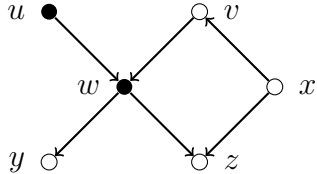

Figure~\ref{fig-ex1} shows a starting position for this game. If Black starts, playing $u$ allows to remove four vertices, i.e. $u,w,y,z$. Then White ends the game by playing $v$ (removing $v$ and $x$). Hence Black wins the game with $4$ points (versus $2$ points for White). On the contrary, if White starts, he can remove all the vertices by playing either $y$ or $z$ as a first move (leaving a unique white vertex that will be removed after his second turn).

Note that this game always ends as removing a single vertex is allowed.

\subsection{Motivation of the game}

We consider two arguers, Black and White, whose objective is to win over a given audience on a number of ideas.~\footnote{\cite{GraMaRu} recently proposed an extension of the de Groot's model of influence based on a similar motivation.}
Following the rules of Aristotelian rhetoric, we assume that the audience and the arguers are endowed with a background of common opinions that can be used by the arguers as premises of their argumentation schemes.

To model the fact that Black and White hold opposite opinions on all subjects, we suppose that Black has a positive opinion on the set $N$ of ideas, and that White has a positive opinion on the negation of each idea in $N$. 
Moreover, we assume that the set $N$ is divided into two disjoint sets $B$ and $W$, representing the ideas that Black and White (respectively) can win through an argument that uses as premises elements of the common background of opinions. 

Each idea is linked to other ideas that can be inferred from it. Correspondingly, once an arguer has won her audience on one idea, she can use it as premise in a new argumentative scheme aimed at winning the audience on some other idea that can be inferred from it. 

This situation can be represented by a directed graph $G = (B \cup W, A)$ where the set of vertices represents the set of ideas in $N$, and an arc $(x, y)\in A$ between two ideas $x$ and $y$ represents the existence of an argument to convince the audience of the veracity of $y$ once $x$ is
accepted as a premise.~\footnote{Our representation of an argumentation system should not be confused with the one proposed by Dung \cite{Dung}. Here, the oriented arcs represent the existence of an argumentative scheme allowing to pass from one idea to another, while in Dung's framework nodes represent ideas (or arguments) and arrows represent attacks between arguments.}

We further assume that whenever there exists an argument that associates to a premise $x$ a conclusion $y$, there exists a ``contrapositive'' argument that associates the negation of $x$ to the negation of $y$.~\footnote{It is straightforward to observe that if the arc from $x$ to $y$ represents a logical implication, the ``contrapositive'' argument is precisely the classical contrapositive proposition.}

To illustrate, suppose that Black aims at convincing the audience on the following ideas:

\begin{tabular}{lcl}
& & \\
$u$: & & We must reduce inequality. \\ 
$v$: & & We must increase the number of officials. \\
$w$: & & We must increase the provision of public goods. \\
$x$: & & We must immediately offer a well-paying job to everyone. \\
$y$: & & We will have to increase tax. \\
$z$: & & Trade needs to be regulated by the government.\\
& &
\end{tabular}


Next, suppose that, from the background of common opinion, Black can convince the audience that ideas $u$ and $w$ are true, while White can convince it that the remaining ideas are false. It follows that $u$, $w\in B$ while $v, x, y, z \in W$.

This situation is depicted in Figure \ref{fig-ex1}, where black edges represent the ideas in $B$, white edges the ideas in $W$, and arrows represent the available argumentation schemes. In particular, Black cannot use $x$ as a premise in her argumentation, but if she convinces the audience that $w$ is true, then she can convince it that $y$ is also true. Similarly, White can convince the audience of idea $\neg w$ only after convinving it of $\neg y$ or $\neg z$.

If Black is the first speaker, she has to pick an idea that she can win by using elements of the background set of common opinions as premises, e.g. $u$. Then, she is able to push through the ideas that can be deduced from $u$, i.e. $w$, $z$ and $y$. 

In the second stage of the game, White can argumentate in order to convince the audience that $v$ is false. By deduction, it will then follow that the audience will also consider $x$ to be false. 

It is clear that the audience's final opinion concerning the ideas in $N$ strongly depends on what player speaks first. 

Another important feature of the example is that, given the background of common opinions, the audience may be convinced about $w \in B$ but against $y \in W$ even though $(w, y) \in A$. For instance, a population of agents may well agree on an increase of public goods ($w$), and yet disagree about paying more taxes ($y$) because, say, this population does not know the tax rate that is needed for the provision of public goods. Such situation may occur because the audience does not know the arguments that can be employed by the arguers beforehand. Alternatively, it can be compared to phenomena such as the doctrinal paradoxes that are common in the literature on judgement aggregation.~\footnote{See for instance \cite{ListPuppe}.}

\subsection{Scoring-play combinatorial game theory}

According to the definition of the ruleset, our game belongs to the universe of scoring-play combinatorial game theory. Roughly speaking, it consists of two-player games with perfect information, no chance and where the moves are alternated. Points are accumulated during the game, and the player having the highest score is the winner.
This class of games has been briefly introduced in the 1950s by Milnor and Hanner in two papers \cite{hanner,milnor}. Their work has been considered as the source of Combinatorial Game Theory - where there is no score and the winner is the player making the last move - that emerges in the following decades (see \cite{Siegel} for a survey book). In parallel, scoring games have been widely considered in the context of economic game theory but there were actually few studies pursuing the work of Milnor and Hanner until the 2010s. Recent years have seen renewed interest for scoring combinatorial games, with attempts to build a general theory as Conway did for combinatorial games. In particular, the PhD theses of Ettinger \cite{ettinger} and Stewart \cite{stewart}, as well as the recent results of Larsson, Nowakowski and Santos~\cite{guaranteed,waiting}, extend the results of Milnor and Hanner for more general classes of scoring games. In addition, new scoring games have been recently considered on a graph structure (graph grabbing game ~\cite{grabbing}, median graph game ~\cite{median}, graph occupation game ~\cite{occupation}, among others, see a larger list in ~\cite{Larsson}). But such games played on graphs have been solved with ad hoc techniques, and apart from the existing scoring game framework. One of our objectives is to study the game \Jeu in correlation with the general theory, and in particular according to the framework defined by Milnor and Hanner.

\subsection{Our results}

In the current paper, we obtain  some properties of  \Jeu.
First, we prove that   \Jeu  is a nonzugzwang game, that means that no player has interest to pass at any  step of the game  (Section \ref{Milnor}). This result seems very natural  but is far from trivial and is very helpful in the following analyses, since it follows that  \Jeu  belongs to Milnor's universe  \cite{milnor}.

As the global study of the game in its whole generality seems to be very difficult, we decided to
restrict it to classes of  directed graphs with some strong properties of symmetry and low degree.
In particular, we introduced the class of \emph{segments} formed by paths of alternated black and white vertices,  with all edges directed from black to white vertices. Even in such a restricted framework, this study remains very complex.
Our main result (Theorem \ref{prop:simple}) claims that for any segment as defined above, the first player always gets a strictly better score that the second one, and the difference between two scores is at most 5.  In other words, we can produce a nearly optimal strategy, in terms of score, for each player for \Jeu played on any segment.

We finish by highlighting and quantifying an interesting property of the game \Jeu. In terms of score, the arrangement of the vertices is more important than the number of vertices of each color. To illustrate this, we exhibit instances of  \Jeu having an arbitrarily  low proportion of black vertices, and where the final proportion of vertices gained by Black is arbitrarily close to 1.

\section{Formalizations}

\subsection{Definitions and notations}

By definition, \Jeu clearly belongs to the universe of scoring-play games. It is partizan, as both players do not always have the same moves because of the two colors. To be consistent with the recent literature about scoring games \cite{Larsson}, we will consider in the rest of the paper the following definition of the game, where the two players are called {\em Left} and {\em Right}, and the vertices are labeled with $L$ and $R$ (instead of Black and White, respectively). \\

The game \Jeu is played on a directed graph $G=(L\cup R, A)$ where the vertex set is partitioned into two sets $L$ and $R$. The vertex set  $L \cup R$ of $G$ will be often denoted by $V(G)$, and to avoid ambiguities about the considered  directed graph,  the sets  $L$, $R$  and $A$ will sometimes be respectively denoted by $L(G), R(G)$ and $A(G)$.\\

For any set $S \subseteq V(G)$,  $G_S $
  denotes the  directed graph $(S, A(G)  \cap S^2)$, and  $G \setminus S$ denotes the directed graph $G_{V(G) \setminus S}$.

  For any $x \in V(G) $, $\mbox{Pred}(G, x)$ denotes the set of vertices that can reach $x$ through a directed path in $G$ (including $x$),
and $\mbox{Succ}(G, x)$ denotes the set of vertices that are reachable from $x$ through a directed path in $G$ (including $x$).\\

The rules of \Jeu are the following. Two players, Left and Right, play alternately by choosing a vertex of their own label:

\begin{itemize}
  \item if it is Left's turn, and $x \in L$ is chosen, then remove $\mbox{Succ}(G, x)$;
  \item if it is Right's turn, and $y \in R$ is chosen, remove  $\mbox{Pred}(G, y)$.
\end{itemize}

The game ends when all the vertices have been removed. The score of each player is the number of vertices he has removed. \\

In the resolution of any scoring game, we are interested in the maximum score gained by each player when both of them play optimally. Such a score will be denoted $s^L_{i}(G)$ and $s^R_{i}(G)$ with $i\in\{1,2\}$, respectively for Left and Right. The value $i$ indicates if the player starts the game or not. For example, $s^L_2(G)$ corresponds to the maximum score gained by Left when Right starts the game. \\

By definition of the ruleset, \Jeu belongs to the class of constant-sum games as we have:
$$
s^R_2(G)+s^L_1(G)=|V(G)|
\quad
\mbox{ and }
\quad
s^L_2(G)+s^R_1(G)=|V(G)|.
$$

\begin{remark}
Note that the scores are defined as the best values reached by each player when considering all the possible sequences of moves. Therefore, the following recursive characterization can be also considered to define the scores of the players:\\

$$s^L_1(G) = s^L_2(G) = 0 \mbox{ when  } |L(G)|=0.
$$
Otherwise,
$$
s^L_1(G) =
\max_{x\in L}\bigl\{|\mbox{Succ}(G,x)|+s^L_2(G\setminus \mbox{Succ}(G,x))\bigr\}.
$$
and
$$
s^L_2(G)
= \min_{y\in R}\bigl\{s^L_1(G\setminus \mbox{Pred}(G,y))\bigr\}.
$$
\end{remark}

As defined in the paper of Milnor~\cite{milnor} or more recently in the survey of Larsson et al. \cite{Larsson}, a more natural way to consider the outcome of a scoring game is through the relative score, i.e. the difference between the scores of Left and Right. As several definitions have been considered in the literature, we choose to use the most recent one defined in \cite{Larsson}. This difference may have two values, called {\em Left score} and {\em Right score}, according to whether Left or Right starts the game.

\begin{definition}
The Left-score and the Right-score of a game $G$, respectively denoted by $Ls(G)$ and $Rs(G)$, are defined as:
$$
Ls(G)=s^L_1(G)-s^R_2(G),
$$
and
$$
Rs(G)=s^L_2(G)-s^R_1(G).
$$
\end{definition}

In other words, games with $Ls(G)>0$ (resp. $Rs(G)>0$) are winning for Left when Left starts (resp. when Right starts) the game. On the contrary, games having a negative relative score are winning for Right.

\begin{ex}
Consider the  directed graph of Figure~\ref{fig-ex1}, where $L$ is the set of black vertices, and $R$ corresponds to the white ones. One can show that the moves described in the introduction are optimal, leading to the following equalities:
\begin{itemize}
\item $s^L_1(G) = 4$ and  $s^R_2(G) = 2$
\item $s^R_1(G) = 6$ and  $s^L_2(G) = 0$
\item $Ls(G)=2$ and $Rs(G)=-6$
\end{itemize}
Since $Ls(G)>0$ and $Rs(G)<0$, each player wants to be the first to move in order to win the game.
\end{ex}

\begin{remark}
As for the absolute scores, when $L(G)$ and $R(G)$ are not empty, we have the following recurrence equations:

$$
Ls(G)=  \max_{x \in L}\bigl \{ \mbox{Succ}(G, x) + Rs( G \setminus \mbox{Succ}(G, x)) \bigr\},
$$
and
$$
Rs(G)=  \min_{y \in R} \bigl \{ -\mbox{Pred}(G, y) + Ls( G \setminus \mbox{Pred}(G, y)) \bigr\},
$$
with the initial condition
$$Ls(G) = - \vert V(G) \vert \  (\mbox{respectively } Rs(G)  = \vert V(G) \vert )$$
when $L(G)$ (respectively $R(G)$) is  empty.
\end{remark}

In the game \Jeu, a first useful property about $Ls(G)$ and $Rs(G)$ is the following:
\begin{lemma}
\label{parity}
In \Jeu, for every game $G$, the quantities $Ls(G)$, $Rs(G)$ and $|V(G)|$ have the same parity.
\end{lemma}
\begin{preuve}
By definition of the scores, we have  $$Ls(G) =   s^L_1(G)-s^R_2(G)  =  2s^L_1(G)  -  |V(G)|$$
and
 $$Rs(G) =   s^L_2(G)-s^R_1(G)  =  2s^L_2(G)  -  |V(G)|.$$
\end{preuve}

Note that in combinatorial game theory, it is common to mix-up a game with a position. Hence $G$ will be used here to indifferently define a game position and its corresponding directed graph.

\subsection{Relevant  directed graphs} \label{relevant}

It may happen that some vertices will be gained by a fixed player, whatever the strategies of both players. These vertices are not relevant from a strategic point of view, hence we should not take them into account. It is the main  idea of this section.\\

More precisely, let  $x \in L$ such that  $\mbox{Succ}(G, x) \subseteq L$.
According to the rules of the game, one can remark  that $x$ and the whole set $\mbox{Succ}(G, x)$ will be removed by Left at some step of the game regardless of Right's future moves. We say that such a vertex $x \in L$ is \emph{ forced} in $G$.


Thus, from a strategic point of view, it is natural to restrict ourselves to directed graphs with no forced vertices\footnote{Actually, there could be a unique difference: if the forced vertices are not immediately removed from the directed graph,  they could be used by a player to pass his turn.  But we will see later (Corollary \ref{cor:nopass}) that passing  is of no interest in this game.}.

\begin{definition}
For each directed graph $ G= (L \cup R,  A)$, let
$$\mbox{Forc}_L(G) =  \bigl\{ x \in L:  \mbox{Succ}(G, x)    \subseteq L\bigr\},  \quad \mbox{Forc}_R(G) =  \bigl\{ y \in R:  Pred(G, y)   \subseteq R\bigr \}$$
 be the subsets of vertices that Left (respectively Right) will remove regardless of the other player's moves, and let
 $$\mbox{Forc}(G) = \mbox{Forc}_L(G) \cup \mbox{Forc}_R(G) $$
be the union of these two subsets. The class of {\bf relevant  directed graphs} for the game is the subclass $\mathcal{S}$ of directed graphs $G = (L \cup R,  A)$ such that $$\mbox{Forc}(G) = \emptyset. $$
\end{definition}

\begin{figure}[h]
			\centering
			\begin{tikzpicture}
			\node[vertex] (0) at (0,2) {};
			\node[vertex] (1) at (2,2) {};
			\node[vertexW] (2) at (1,1) {};
			\node[vertex] (3) at (3,1) {};
			\node[vertexW] (4) at (0,0) {};
			\node[vertexW] (5) at (2,0) {};
			\draw[edge](0)--(2);
			\draw[edge](1)--(2);
			\draw[edge](2)--(4);
			\draw[edge](2)--(5);
			\draw[edge](3)--(5);
                        \draw[edge](0)--(4);
                        \draw[edge](1)--(3);
            \draw(0) node[left=5] {$x$};
             \draw(1) node[right=5] {$x'$};
			
			\end{tikzpicture}
			
	\caption{The  directed graph $G$ is relevant but after playing $x\in L$, the new directed graph $G\setminus \mbox{Succ}(G,x)$ is not anymore relevant since $x'$ is now forced.}
	\label{fig-ex2}
	\end{figure}
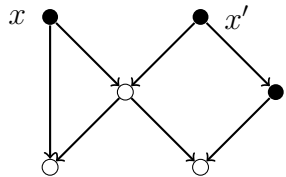

For example, for the directed graph of Figure~\ref{fig-ex1}, the vertices $v$ and $x$ are forced, as player White will always get them at the end.
Remark that,  for each directed graph $G$, we have
$ G \setminus \mbox{Forc} (G) \in \mathcal{S}.$
Also note that $G \in \mathcal{S}$ does not imply that
$G \setminus  \mbox{Succ}(G, x) \in \mathcal{S}$, since  $G \setminus  \mbox{Succ}(G, x)$  can contain some forced vertices, see Figure~\ref{fig-ex2}.
Nevertheless, it can be decided to immediately attribute these forced vertices (which are necessarily in $L$) to Left. This allows to have a new relevant directed graph:
$$(G \setminus  \mbox{Succ}(G, x)) \setminus  \mbox{Forc} (G \setminus  \mbox{Succ}(G, x)), $$
 and, therefore,
 to have recursive formulas  for relevant directed graphs. This is done using  the following definitions.

\begin{definition}
Let  $G \in  \mathcal{S}$ and  $x \in L$. The set $\mbox{Rmv}(G, x)$ is  the set defined as:
$$ \mbox{Rmv}(G, x) =  \mbox{Succ}(G, x) \cup \mbox{Forc} (G \setminus  \mbox{Succ}(G, x)).$$
The  directed graph $G_x$ defined as:
$$G_x =  G \setminus   \mbox{Rmv}(G, x)$$
is an element of  $\mathcal{S}$. This  directed graph  is called the {\bf  resulting relevant directed subgraph induced by $x$}.
In a similar way, for $y \in R$, we set
$$  \mbox{Rmv}(G, y) =  \mbox{Pred}(G, y)  \cup \mbox{Forc} (G \setminus  \mbox{Pred}(G, y)).$$
and
$$G_y =  G \setminus   \mbox{Rmv}(G, y).$$
\end{definition}

\begin{remark} \label{rmk:rmv} We can alternatively define $ \mbox{Rmv}(G, x) $ and $\mbox{Rmv}(G, y) $ by

$$ \mbox{Rmv}(G, x) = \bigl\{z \in V(G) \ |\   \mbox{Succ}(G, z)  \subseteq L\, \cup \, \mbox{Succ}(G, x) \bigr\}$$ \mbox{ and }$$ \mbox{Rmv}(G, y) = \bigl\{z' \in V'(G)\ |\ \mbox{Pred}(G, z' )  \subseteq R \,\cup \, \mbox{Pred}(G, y) \bigr\}. $$
\end{remark}

 \textbf{Recursive formulas. }
In order not to leave the class $\mathcal{S}$, we adapt the recursive formulas of the relative scores.
Precisely, for each nonempty $G \in \mathcal{S}$, since after a move of a player, the forced vertices are necessarily for this player, we have:
 $$s^L_1(G)=\max_{x\in L} \bigl \{ |\mbox{Rmv}(G,x)|+s^L_2(G_x) \bigr\} \quad \mbox{and}    \quad s^L_2(G)=\min_{y \in R} \bigl\{s^L_1(G_y) \bigr\},$$

$$s^R_1(G)=\max_{y\in R} \{ |\mbox{Rmv}(G,y)|+s^R_2(G_y)\} \quad \mbox{and}    \quad s^R_2(G)=\min_{x \in L} \bigl\{s^R_1(G_x) \bigr\},$$

$$
Ls(G)=  \max_{x \in L}\bigl\{ \vert \mbox{Rmv}(G, x) \vert  -
Rs(G_x)\bigl\} \quad \mbox{and} \quad Rs(G)=  \min_{y \in R} \bigl\{- \vert \mbox{Rmv}(G, y) \vert   +
Ls(G_y)\bigl\}.
$$

with, for initialization, the convention that for the empty graph all these values are null.

\section{\Jeu belongs to Milnor's universe of scoring games}
\label{Milnor}

In the literature, different universes of scoring games have been defined in the last decades, according to structural properties of the game. In particular, Milnor and Hanner restricted their study to {\em nonzugzwang dicot} games:
\begin{itemize}
\item A game has no zugzwang if each player always prefers moving rather than missing his turn. In other words, it satisfies $Ls(G)\geq Rs(G)$ for every subposition $G$ of the game.
\item A game is a dicot if, for every subposition of the game, a player can move if and only if his opponent also can.
\end{itemize}

In Milnor's paper, the {\em incentive} of a game $I(G)$ is defined as the value $Ls(G)-Rs(G)$. It corresponds to the interest that Left has to move rather than passing. Hence games in Milnor's universe are also called dicot games with nonnegative incentive.

It turns out that this family of scoring combinatorial games has a nice mathematical behavior that induces a group structure, following the case of normal-play combinatorial games. For this group, the binary operation $(+)$ is called the {\em sum of games}. Given two games $G$ and $G'$, the game $G+G'$ is the game where each player moves either in $G$ or in $G'$. The games ends when all the moves are exhausted in both $G$ and $G'$.\\

We now claim that the game \Jeu played on each relevant  directed graph $G$, using moves induced by sets $\mbox{Rmv}(G, x)$, can be embedded into Milnor's universe. First, from the definition of relevant  directed graphs, \Jeu is a dicot game.  The key-point is the fact that \Jeu is a nonzugzwang game (Corollary \ref{cor:nopass}).

\subsection{Basic Lemmas}
The proof is based on three lemmas. The first one describes the evolution of the succesors and predecessors of a vertex after a move is played.

\begin{lemma}\label{lem:evolsucc}
 Let $G\in \mathcal{S}$ and let $u,v \in V(G)$ such that $v\in V(G_u)$. Then, $$\mbox{Succ}(G_u,v)=\mbox{Succ}(G,v)\cap V(G_u)$$ and $$\mbox{Pred}(G_u,v)=\mbox{Pred}(G,v)\cap V(G_u).$$ 
  \end{lemma}

\begin{preuve}
Since there is no assumption on $u$ and $v$ about being in $L$ or $R$, the two results are symmetric. Thus we just prove the first equality.
  
 Let $z\in V(G_u)$. If $z$ is a successor of $v$ in $G_u$, it is also a successor of $v$ in $G$. Thus $\mbox{Succ}(G_u,v)\subseteq \mbox{Succ}(G,v)\cap V(G_u)$. Assume now that $z$ is a successor of $v$ in $G$. If it is not a successor of $v$ in $G_u$, it means that on any path from  $v$ to  $z$ in $G$, there is a vertex in $\mbox{Rmv}(G,u)$. Let $z'$ be such a vertex. If $u\in L$, then it means that $\mbox{Succ}(G,z')\subseteq L\cup \mbox{Succ}(G,u)$. In particular, since $z$ is a successor of $z'$,  $\mbox{Succ}(G,z)\subseteq  \mbox{Succ}(G,z')  \subseteq L\cup \mbox{Succ}(G,u)$, so  $z \in \mbox{Rmv}(G,u)$, a contradiction since, by assumption, $z\in G_u$. Otherwise, if $u\in R$, then $\mbox{Pred}(G,z')\subseteq R\cup \mbox{Pred}(G,u)$. In particular, since $v$ is a predecessor of $z'$,  $\mbox{Pred}(G,v)\subseteq \mbox{Pred}(G,z')\subseteq R\,\cup\, \mbox{Pred}(G,u)$, so  $v \in \mbox{Rmv}(G,u)$, a contradiction.

\end{preuve}

The next lemma describes what happens when two moves are commuted.

\begin{lemma}\label{lem:commute} 
For each $G \in  \mathcal{S}$, each $x, x'  \in L$ and each $y \in R$, we have:
\begin{enumerate}
\item   If $  x' \in  \mbox{Rmv}(G, x)$,  then $\mbox{Rmv}(G, x') \subseteq \mbox{Rmv}(G, x)$, and $(G_{x'})_x = G_x$.
\item    If $  x' \notin  \mbox{Rmv}(G, x)$  and  $x \notin  \mbox{Rmv}(G, x')$,  then  $(G_{x'})_x = (G_x)_{x'}$. 
\item    If $  x \notin  \mbox{Rmv}(G, y)$,  then  $\mbox{Rmv}(G_x, y)=\mbox{Rmv}(G, y)$ and $(G_{x})_y = (G_y)_{x}$.
\end{enumerate}
\end{lemma}

Before turning to the proof of Lemma \ref{lem:commute}, it is worth mentioning that, because of the partition  of $V(G)$ in  $L$ and $R$,  the condition   $  x \notin  \mbox{Rmv}(G, y)$ is actually symmetric in $x$ and $y$. Precisely, we have:
$$  \bigl [x \in  \mbox{Rmv}(G, y) \bigr]  \iff \bigl [ x \in  \mbox{Pred}(G, y)\bigr] \iff \bigl [ y \in  \mbox{Succ}(G, x) \bigr]  \iff \bigl [y \in  \mbox{Rmv}(G, x)\bigr].$$

\begin{preuve}
  
\noindent  {\bf 1.}
  Assume that  $ x' \in  { \mbox{Rmv}(G, x)}$ and let $z \in \mbox{Rmv}(G, x')$.   Then, from Remark  \ref {rmk:rmv},  $ \mbox{Succ}(G, z) \subseteq L\, \cup \, \mbox{Succ}(G, x')$.  But, since $ x' \in  \mbox{Rmv}(G, x)$, we  have $\mbox{Succ}(G, x')  \subseteq L\, \cup \, \mbox{Succ}(G, x)$. Thus,  $ \mbox{Succ}(G, z) \subseteq L\, \cup \, \mbox{Succ}(G, x)$ which gives that  $z \in \mbox{Rmv}(G, x)$ and so $\mbox{Rmv}(G, x') \subseteq \mbox{Rmv}(G, x)$, as desired. 
To show that  $(G_{x'})_x = G_x$,  we have to prove that
$ \mbox{Rmv}(G, x) =   \mbox{Rmv}(G_{x'}, x) \cup   \mbox{Rmv}(G, {x'})$,
or equivalently, 
$$ \mbox{Rmv}(G, x) \setminus  \mbox{Rmv}(G, {x'})=   \mbox{Rmv}(G_{x'}, x), $$
where the equivalence follows from the fact that  $\mbox{Rmv}(G_{x'}, x) \cap   \mbox{Rmv}(G, {x'}) = \emptyset$ and $\mbox{Rmv}(G, {x'}) \subseteq \mbox{Rmv}(G, {x})$.
First, let $z \in  \mbox{Rmv}(G_{x'}, x)$. By Lemma \ref{lem:evolsucc}, $\mbox{Succ}(G_{x'}, x ) =  \mbox{Succ}(G, {x}) \cap  V(G_ {x'}) $, and so
$$\mbox{Succ}(G_{x'}, z )  \subseteq   L \cup \mbox{Succ}(G_{x'}, x )  \subseteq    L \cup \mbox{Succ}(G, {x}), $$
 which proves that $z \in  \mbox{Rmv}(G, x) $. Obviously,  $z \in V(G_{x'})$ and so $z \notin  \mbox{Rmv}(G, {x'})$. From this, we finally get that $z \in \mbox{Rmv}(G, x) \setminus  \mbox{Rmv}(G, {x'})$. This proves that $\mbox{Rmv}(G_{x'}, x) \subseteq \mbox{Rmv}(G, x) \setminus  \mbox{Rmv}(G, {x'})$.

Conversely, let $ z \in  \mbox{Rmv}(G, x) \setminus  \mbox{Rmv}(G, {x'})$.
We have:
\begin{eqnarray}
\mbox{Succ}(G_{x'}, z )  &= &    \mbox{Succ}(G, {z}) \cap V(G_{x'} ) \nonumber \\
& \subseteq  & \bigl(L \cup \mbox{Succ}(G, {x}) \bigr) \cap V(G_{x'})  \nonumber \\
& \subseteq  & L \cup \bigl(\mbox{Succ}(G, {x})  \cap V(G_{x'}) \bigr) \nonumber \\
& =  & L \cup \mbox{Succ}(G_{x'}, {x}), \nonumber
\end{eqnarray}
which proves that $z \in  \mbox{Rmv}(G_{x'}, x )$ and completes the proof of Point 1. \\

\

\noindent {\bf 2.}  Assume that $  x' \notin  \mbox{Rmv}(G, x)$  and  $x \notin  \mbox{Rmv}(G, x')$, and let $ \mbox{Rmv}(G, x, x')$ be the set of vertices defined as:
$$ \mbox{Rmv}(G, x, x') =   \bigl\{z\in V(G)\ | \ \mbox{Succ}(G, z) \subseteq L   \cup  \mbox{Succ}(G, x)   \cup  \mbox{Succ}(G, x') \bigr\}.   $$
We will prove that
$$(G_{x'})_x = G \setminus \mbox{Rmv}(G, x, x') $$ which
 is equivalent to prove that:
\begin{equation} \label{eq1}
\mbox{Rmv}(G, x, x')  =  \mbox{Rmv}(G, x') \cup  \mbox{Rmv}(G_{x'}, x)
\end{equation}
If
$ z \in  \mbox{Rmv}(G, x')$,  then $\mbox{Succ}(G, z) \subseteq L\, \cup \, \mbox{Succ}(G, x' )$, which yields that  $z \in  \mbox{Rmv}(G, x, x')$.
If $ z \in  \mbox{Rmv}(G_{x'}, x)$,  then  $\mbox{Succ}(G_{x'}, z) \subseteq L\, \cup \, \mbox{Succ}(G_{x'}, x)$, or equivalently, using Lemma \ref{lem:evolsucc},
 $$\mbox{Succ}(G, z) \cap  V(G_{x'}) \subseteq L\, \cup \,( \mbox{Succ}(G, x) \cap  V(G_{x'})).$$
 This leads to
 $$\mbox{Succ}(G, z) \cap  V(G_{x'})   \subseteq L\, \cup \, \mbox{Succ}(G, x). $$
On the other hand,
 $\mbox{Succ}(G, z) \setminus   V(G_{x'} )=  \mbox{Succ}(G, z) \cap   \mbox{Rmv}(G, x')$.
Because $\mbox{Rmv}(G, x') \subseteq L \cup  \mbox{Succ}(G, x')$, we have  $$\mbox{Succ}(G, z)   \setminus   V(G_{x'})  \subseteq L \cup  \mbox{Succ}(G, x'). $$
It follows that $\mbox{Succ}(G, z) \subseteq L \cup  \mbox{Succ}(G, x) \cup  \mbox{Succ}(G, x'),  $ that is, $z \in
\mbox{Rmv}(G, x, x'). $ We have proved that $\mbox{Rmv}(G, x, x')  \supseteq  \mbox{Rmv}(G, x') \cup  \mbox{Rmv}(G_{x'}, x).$

For the converse inclusion, take $z \in  \mbox{Rmv}(G, x, x')$, that is, $z$ is such that  $\mbox{Succ}(G, z) \subseteq L \cup  \mbox{Succ}(G, x) \cup  \mbox{Succ}(G, x').$
If  $\mbox{Succ}(G, z) \subseteq L \cup    \mbox{Succ}(G, x')$, then $z \in  \mbox{Rmv}(G, x')$.
Otherwise, $z \in  V(G_{x'})$  and
\begin{eqnarray}
\mbox{Succ}(G_{x'}, z) &= & \mbox{Succ}(G, z) \cap V(G_{x'})  \nonumber \\
&\subseteq & \bigl(L \cup  \mbox{Succ}(G, x) \cup \mbox{Succ}(G,x')\bigr)  \cap V(G_{x'}) \nonumber \\
&\subseteq & L   \cup \bigl( \mbox{Succ}(G, x)  \cap V(G_{x'})\bigr)   \nonumber \\
&=& L \cup \mbox{Succ}(G_{x'}, x).  \nonumber
\end{eqnarray}
Thus $z \in  \mbox{Rmv}(G_{x'}, x)$.
  This proves claim (\ref{eq1}). The result $(G_{x'})_x =  (G_{x})_{x'}$ obviously follows. \\

  \
  
\noindent {\bf 3.} Let    $x \in L$ and $y \in R$  such that  $x \notin  \mbox{Rmv}(G, y)$. We claim that
  $$ \mbox{Rmv}(G, x) \cap  \mbox{Rmv}(G, y)  = \emptyset.$$
  Indeed, assume it is not the case and let $z\in  \mbox{Rmv}(G, x) \cap  \mbox{Rmv}(G, y)$. Without loss of generality, assume that $z\in L$. Then since $z\in  \mbox{Rmv}(G, y)$, we must have $z\in \mbox{Pred}(G,y)$. But then $y$ is a successor of $z$, and thus must be in $\mbox{Rmv}(G,x)$ since $z$ is, leading to a contradiction.


  We now prove that $$\mbox{Rmv}(G_y, x) = \mbox{Rmv}(G,x).$$
  Let $z\in V(G_y)$. Because $z\in V(G_y)$ implies that all its successors are also in $V(G_y)$,  by Lemma \ref{lem:evolsucc} we have  $\mbox{Succ}(G_y, z) = \mbox{Succ}(G, z)\cap V(G_y) = \mbox{Succ}(G,z)$. In particular, this is also true for $x$. Thus, we have $\mbox{Succ}(G_y, z) \subseteq L\ \cup \ \mbox{Succ}(G_y,x)$ is equivalent to $\mbox{Succ}(G, z) \subseteq L\ \cup \ \mbox{Succ}(G,x)$,  and thus $\mbox{Rmv}(G_y, x) = \mbox{Rmv}(G,x)$.

Symmetrically, we have $\mbox{Rmv}(G_x, y) =\mbox{Rmv}(G, y)$. And finally, we get
 $$ (G_y)_x =  G \setminus \bigl( \mbox{Rmv}(G, y) \cup \mbox{Rmv}(G, y)\bigr) = (G_x)_y. $$

\end{preuve}

The third lemma below uses the  previous one to obtain a strong result about the consequences of a move.

\begin{lemma}\label{lem:subposition2}
  Let $G=(L\cup R,A)$ be a  directed graph of $\mathcal{S}$. Let $y\in R$.  
   We have
  $$\forall i \in \{1, 2\}, \quad s_i^L(G)\geq s_i^L(G_y).  $$
  Similarly, let $x \in L$. We have
  $$\forall i \in \{1, 2\}, \quad s_i^R(G)\geq s_i^R(G_x).$$
   \end{lemma}

\begin{preuve}
We prove both results by induction on $G$. First, we consider the case $y\in R$.
If $G_y $ is the empty graph, the result is clear since $s^L_1(G_y) = s^L_2(G_y)=0$.
Otherwise, we distinguish two cases.

\noindent {\bf Case 1.} Suppose that $i=1$. Recall that we have:
 $$s^L_1(G)=\max_{x\in L} \bigl\{ |\mbox{Rmv}(G, x)|+s^L_2(G_x)\bigr\}.$$

We denote by $L_{G_y} = L \cap V(G_y)$,  the set of  Left  vertices of the directed subgraph $G_y$.
 We have
  $$s^L_1(G_y)= \max_{x\in L_{G_y} } \bigl \{|\mbox{Rmv}(G_y, x)|+s^L_2((G_y)_x) \bigr \}.$$
From Lemma \ref{lem:commute}, we have
 $ (G_y)_x =  (G_x)_y $, thus
we obtain
 $$s^L_1(G_y)  = \max_{x\in L_{G_y}} \bigl\{|\mbox{Rmv}(G_y,x)|+s^L_2((G_x)_y)\bigr\}.$$
   By the induction hypothesis, we have $s^L_2((G_x)_y)\leq s^L_2(G_x)$. Moreover, $\mbox{Rmv}(G_y,x)=\mbox{Rmv}(G,x)$ since $x\in G_y$. We finally obtain:
$$s^L_1(G_y)   =  \max_{x\in L_{G_y}} \{|\mbox{Rmv}(G_y,x)|+s^L_2((G_x)_y)\}
  \leq \max_{x\in L} \{|\mbox{Rmv}(G,x)|+s^L_2(G_x)\} = s^L_1(G).$$
\noindent {\bf Case 2.} Suppose that $i=2$. We have:
$$s^L_2(G)=\min_{y'\in R} \bigl\{s^L_1(G_{y'})\bigr\}$$
 Let $y^*\in R$ the vertex that realizes this minimum. We have $s^L_2(G)=s^L_1(G_{y^*})$.

We distinguish three subcases.

\noindent {\bf Case 2-1.}
$y^*\notin \mbox{Rmv}(G,y)$ and $y\notin \mbox{Rmv}(G,y^*)$.
From Lemma \ref{lem:commute}, we have    $(G_{y})_{y^*} = (G_{y^*})_y.  $
Thus, we have:
$$ s^L_1((G_{y^*})_{y}) = s^L_1((G_{y})_{y^*})= s^L_1(G_y \setminus \mbox{Rmv}(G_y, y^*))\geq    \min_{z \in R_{G_y}} \bigl\{s^L_1(G_y \setminus \mbox{Rmv}(G_y ,z))\bigr\} = s^L_2(G_y),  $$
where  $R_{G_y}$ denotes the set of $R$ vertices  of $G_y$, and where the inequality follows from the fact that $y^*$ is not necessarily optimal in $G_y$.
On the other hand, by the induction hypothesis,  we have
$$s^L_1(G_{y^*}) \geq s^L_1((G_{y^*})_y). $$

Thus, by combining both inequalities,  we finally get
 $$s^L_2(G) = s^L_1(G_{y^*}) \geq s^L_1((G_{y^*})_y)\geq s^L_2(G_y),  $$ as desired.

\noindent {\bf Case 2-2.}  $y^*\notin \mbox{Rmv}(G,y)$ but $y\in  \mbox{Rmv}(G, y^*)$. Then, from Lemma \ref{lem:commute},  $G_{y^*}=(G_y)_{y^*}$. It follows that
   $$s^L_2(G) = s^L_1(G_{y^*})=s^L_1((G_y)_{y^*})\geq s^L_2(G_y).$$

\noindent {\bf Case 2-3.}   $y^*\in \mbox{Rmv}(G, y)$. Then, from Lemma \ref{lem:commute},  $G_{y}=( G_{y^*})_y.$ Thus, using the definition of $s^L_2(G)$ and  the induction hypothesis on $G_{y^*}$, we obtain that
 $$ s^L_2(G) = s^L_1(G_{y^*}) \geq s^L_1((G_{y^*})_y)  =  s^L_1(G_y).  $$
Once again,  by induction, for each $z \in R_{G_y} $, $s^L_1(G_y) \geq s^L_1((G_{y})_{z})$. Hence,
   $$s^L_1(G_y)  \geq \min_{z \in R_{G_y}} \bigl\{ s^L_1((G_{y})_{z}) \bigr\} \geq s^L_2(G_y). $$%
 Combining the inequalities above  allows to conclude that  $s^L_2(G) \geq s^L_2(G_y) $.

 This completes the proof for $L$.
Using similar arguments, the proof follows for $R$.
\end{preuve}

\subsection{Consequences}

 Lemma \ref{lem:subposition2} has some useful  corollaries.
The first corollary indicates that    \Jeu is a nonzugzwang game.

\begin{corollary}[nonzugzwang] \label{cor:nopass}
  Let $G=(L\cup R,A)$ be a  directed graph of   $\mathcal{S}$. For each $K\in \{R,L\}$, $s^1_K(G)\geq s^2_K(G)$.
\end{corollary}

\begin{preuve}
By formulas of   Section \ref{relevant},
$$s^L_2(G)=\min_{y\in R} \bigl\{s^L_1(G_y))\bigr\}.$$
By Lemma \ref{lem:subposition2},
 $$\forall y \in R, \quad s^L_1(G_y)\leq s^L_1(G),$$
from which we conclude that $s^L_1(G) \geq s^L_2(G)$.
 The proof that $s^R_1(G)\geq s^R_2(G)$ is similar.
\end{preuve}

The next corollary indicates that for $y', y \in R$ and $y' \in \mbox{Rmv}(G,y)$,
playing $y'$ cannot be better for Right than playing $y$.

\begin{corollary} [included moves]
\label{cor:included}
Consider any directed graph $G=(R\cup L,A)$ of  $\mathcal{S}$.

Let $y, y' \in R$ be such that  $y'  \in \mbox{Rmv}(G,y)$. Then, 
  $$s^L_1(G_{y'} )\geq s^L_1(G_y).$$

 Let $x, x' \in L$ be such that  $x'  \in \mbox{Rmv}(G,x)$. Then,
  $$s^R_1(G_{x'})\geq s^R_1(G_x).$$
\end{corollary}

\begin{preuve}
By Lemma  \ref{lem:subposition2},  we have
$$s^L_1(G_{y'})\geq s^L_1((G_{y'})_y) .$$
Since $y'  \in \mbox{Rmv}(G, y)$, we have
$$(G_{y'})_y = G_y, $$
 from which we conclude that $s^L_1(G_{y'} )\geq s^L_1(G_y)$.
 The proof of the second claim is similar.
 \end{preuve}

Therefore, as the game   \Jeu played on relevant directed graphs is a nonzugzwang dicot, it satisfies the properties of Milnor's universe. In particular,
Milnor provides interesting properties about the sum of two games of this universe.

\begin{corollary}[Milnor~\cite{milnor}]
\label{cor:milnor}
For two nonzugzwang dicot games $G$ and $G'$, we have
$$
\left.
\begin{array}{r}
Ls(G)+Rs(G')   \\
Rs(G)+Ls(G')
\end{array} \right\} \leq Ls(G+G') \leq Ls(G)+Ls(G')
$$
\end{corollary}

Note that a symmetric result is also available for $Rs(G+G')$.

 \begin{preuve} (sketch)
 The first player Left can use the strategy consisting  in first   playing optimally in $G$, and then playing optimally for  the same game where Right has just played.  The  nonzugzwang hypothesis is necessary for the case when  Left  cannot reply in the same game, say $G$, since this game is over.  In this case, the second player Left plays (optimally) in the remaining game $G'$, and by this way, plays two consecutive moves in $G'$.  This does not decrease its score in $G'$, since $G'$ is a nonzugzwang game.  This proves the lower bounds.

On the other hand, the second player Right can reply to each move of Left by the best strategy in the  game where Left just played, $G$ or $G'$.  As for the previous inequality, the  nonzugzwang hypothesis is necessary for the case when Right cannot reply in the same game.  This proves that $Ls(G+G') \leq Ls(G)+Ls(G')$.
\end{preuve}

 The interest of Corollary~\ref{cor:milnor} is natural in the case of \Jeu, since each time a move disconnects the  directed graph, the resulting game corresponds to a sum of smaller games. It will be applied on particular instances of the game.\\

 On the other hand, the  classical result  \cite{milnor}  about the the addition of two opposite  games  can be interpreted as follows  for \Jeu.

\begin{corollary}[Milnor~\cite{milnor}]
\label{cor:G-G}
Let $H$ be a relevant  directed graph,  and $-H$ be the directed graph such that $L$-vertices (respectively $R$-vertices) of $-H$ are $R$-vertices (respectively $L$-vertices) of  $H$ and $(x, y)$ is a directed arc  of $-H$ when  $(y, x)$ is an directed arcs  of $H$.
We have
$$
 Ls(H+(-H)) = Rs(H+(-H))= 0.
$$
\end{corollary}

\begin{preuve} [(sketch)]  The fact that  $Ls(H+(-H))  \leq  0$ can be shown by using a copy strategy: when Left chooses a position in  $H$ (respectively in $-H$), then Right replies by choosing the same position in  $-H$   (respectively in $H$).
Regarding the other inequality, by symmetry we have  $s_1^L(H+(-H)) =  s_1^R(H+(-H))$ and $s_2^L(H+(-H)) =  s_2^R(H+(-H))$. Furthermore, by Corollary \ref{cor:nopass},  we have $s_1^L(H+(-H)) \geq s_2^L(H+(-H)) =  s_2^R(H+(-H))$. It follows that  $ Ls(H+(-H)) = s_1^L(H+(-H)) -  s_2^R(H+(-H)) \geq 0$. The proof is similar for $Rs(H+(-H))$.
\end{preuve}

Note that in combinatorial game theory, the game $-H$ is generally called the {\em negative} of $H$ and is defined as the game $H$ where the roles of Left and Right are interchanged. In particular, we have the following interesting result derived from this definition:
$$
Ls(-H)=-Rs(H)
$$

\section{\Jeu played on unions of alternated directed paths or cycles}

As for many games, it seems very complex to study the game \Jeu in general. Focusing on particular instances is thus natural. First, one can remark that the structure of the game is correlated to strongly connected components of the directed graph: when  a vertex $x$ is chosen by Left, then all vertices of the strongly connected component of $x$  are taken by Left,  the other  strongly connected components containing a successor of  $x$ are taken by Left, while all vertices of the other strongly components are not taken and remain in the current position. Thus, we choose to focus on acyclic graphs.

In  a relevant acyclic  directed graph, all sources are elements of $L$ and all sinks are elements of $R$. Moreover,  from  Corollary  \ref{cor:included},   all  strategic optimal  moves are played on source or sink vertices. After a move at a source $x$, all successors of $x$, and possibly some other vertices in $L$, are removed. The resulting directed graph may be highly complex.
To reduce this complexity, we decide to only consider instances where there are only sources and sinks, i.e. bipartite  directed graphs.
Finally,  we added some regularity hypotheses to introduce the classes of  directed graphs below that correspond to alternated directed paths and cycles.

\begin{definition}
  A finite  directed graph  $G = (L \cup R, A)$ is called a \emph{segment}  if the following conditions hold.
\begin{itemize}
\item   the set $V(G)$ is a set of consecutive  integers (i.e. $(i, j \in V(G) \land  i \leq k \leq  j) \implies   (k \in V(G)) $), such that $\vert V(G) \vert  \geq 2$;
\item  $L$ is   the  set of even integers of $V(G) $, and $R$  is   the  set of odd integers of $V(G) $;
\item if $2i$ and $2i+1$ are both vertices of $G$, then $(2i, 2i+1) \in A$;
\item if $2i$ and $2i-1$ are both vertices of $G$, then $(2i, 2i-1) \in A$;
\end{itemize}

The  lowest and largest integers of  a segment are called the \emph{endpoints}
 of the segment.
The class of segments is denoted by $\mathcal{C}^{segment}$.

A finite  directed graph is  an element of the class $\mathcal{C}$ if it is a finite sum of segments. The endpoints of a directed  graph of $\mathcal{C}$ are the endpoints of the  segments which compose it.

A finite  directed graph  $G = (L \cup R, A)$ is an element of the class $\mathcal{C}^{cycle}$
if it is obtained from a segment with an even number of vertices by adding the  arc $(x, y)$, where $x$  is  the endpoint of the segment which belongs to $L$ and $y$  is  the endpoint of the segment which belongs to $R$
\end{definition}

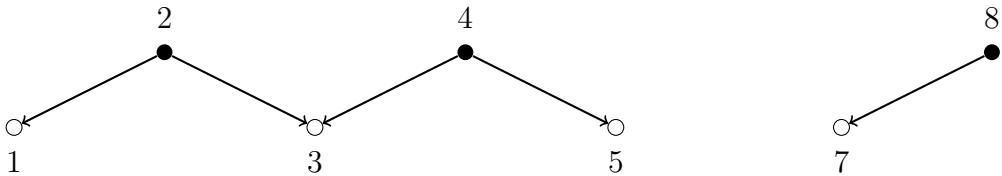
\begin{figure}[h]
			\centering
			\begin{tikzpicture}
            \node[vertexW] (1) at (0,0) {};
			\node[vertex] (2) at (2,1) {};
			\node[vertexW] (3) at (4,0) {};
			\node[vertex] (4) at (6,1) {};
			\node[vertexW] (5) at (8,0) {};
			\node[vertexW] (7) at (11,0) {};
			\node[vertex] (8) at (13,1) {};
			\draw[edge](2)--(1);
			\draw[edge](2)--(3);
			\draw[edge](4)--(3);
			\draw[edge](4)--(5);
			\draw[edge](8)--(7);
            \draw(1) node[below=5] {$1$};
             \draw(2) node[above=5] {$2$};
              \draw(3) node[below=5] {$3$};
               \draw(4) node[above=5] {$4$};
                \draw(5) node[below=5] {$5$};
                 \draw(7) node[below=5] {$7$};
                 \draw(8) node[above=5] {$8$};
			
			\end{tikzpicture}
			
	\caption{The sum of the two segments $[1,5]$ and $[7,8]$}
	\label{fig:segments}
	\end{figure}

Notice that $\mathcal{C} \subseteq \mathcal{S}$ and $\mathcal{C}^{cycle}  \subseteq \mathcal{S}$. Figure~\ref{fig:segments} illustrates the class $\mathcal{C}$.

\subsection{Values reached for sums of segments}

\textbf{Classification of moves in $\mathcal{C} $ } For each integer $k$, a {\bf $k$-move} is a position such that a player can remove $k$ vertices from the directed  graph,  with the convention that forced  vertices created by the move are immediately removed.
In a  directed graph $G \in \mathcal{C}$, there are four possibilities:
\begin{itemize}
\item a 2-move. Such a move is initiated on a endpoint of a segment whose  size is different from 3. Indeed, if the  size of the segment were 3, a forced vertex would be  created. On Figure~\ref{fig:segments}, playing the vertex $1$, $5$, $7$ or $8$ is a $2$-move.
\item a 3-move. It appears when one the following situations holds:
\begin{itemize}
\item (a) the move is initiated on a segment of size 3;
\item (b) the move is initiated on a segment of size at least 5 from a vertex $i$ such that either $i-1$ or $i+1$ is an endpoint of the segment, but not both (e.g. playing $2$ or $4$ on Figure~\ref{fig:segments}).
\item (c) the move is initiated from a vertex $i$ such that $ i-3, i-2, i-1, i+1, i+2, i+3$ are vertices of $G$;
\end{itemize}
\item  a 4-move. Such a move is initiated on  a segment of size at least 4, from a vertex $i$ such that  $ i-1, i+1$ are vertices of $G$, and either $i-2$ or $i+2$ is an endpoint of $G$, but not both.
\item a 5-move. The move is initiated from a vertex $i$ such that $i-2, i-1, i+1, i+2$ are vertices of $G$ and  $i-2$ or $i+2$ both are endpoints of $G$. It is the case of playing the vertex $3$ on Figure~\ref{fig:segments}. Such moves can only appear on segments of size $5$.
\end{itemize}

There is no other kind of moves when $G \in \mathcal{C}$. Except in case (c) in a 3-move, an endpoint is always removed  by the active player and thus the number of segments does not increase.  A 3-move of type (c) is called a {\em cutting} move while all other moves are called {\em border} moves.

\begin{remark}\label{rem:2move}
  Note that in an optimal strategy, $2$-moves are only played on segments of size $2$. Indeed, if a $2$-move is played on a larger segment, then there is always in this segment a $k$-move, with $k\geq 3$ that contained the $2$-move and thus, by Corollary \ref{cor:included}, will lead to a better score. 
\end{remark}

%
%

The next proposition limits the value of Left and Right score when there is at most one segment with an odd number of vertices. Before stating this result we need a definition.

\begin{definition}  We distinguish three subsets of segments.
\begin{enumerate}

\item $ \mathcal{C}_=$ stands for the subset of $ \mathcal{C}$ formed by  directed graphs $G$ whose all segments have an even number of vertices.
\item $ \mathcal{C}_+$ represents the subset of $ \mathcal{C}$ formed by directed  graphs $G$ whose all segments have an even number of vertices except one,  and $\vert L \vert -  \vert R \vert = 1$.
\item $ \mathcal{C}_-$ stands for the subset of $ \mathcal{C}$ formed by directed  graphs $G$ whose all segments have an even number of vertices except one,  and $\vert L \vert -  \vert R \vert = -1$.
\end{enumerate}
\end{definition}

\begin{proposition}
 \label{prop:segments}
We have the following results:
\begin{enumerate}
\item   For each $G \in \mathcal{C}_=$, we have $Ls(G) \in \{0, 2, 4\}$ and $Rs(G) \in \{0, -2, -4\}$.
\item   For each $G \in \mathcal{C}_+$, we have   $Ls(G) \in \{1, 3, 5\}$ and $Rs(G) \in \{1, -1, -3\}$.
\item   For each $G \in \mathcal{C}_-$, we have   $Ls(G) \in \{-1, 1, 3\}$  and $Rs(G) \in \{-1, -3, -5\}$.
\end{enumerate}
\end{proposition}


\begin{preuve}  By induction on the number $n = \vert L \vert + \vert R \vert $ of vertices of $G$.
The main principle of the induction is to preserve instances of $ \mathcal{C}$ containing at most one segment with an odd number of vertices.

For initialization, an  easy exhaustive analysis  shows that Proposition \ref{prop:segments} is true for $ n \leq 5$.
Next, take a  directed graph $G$ of  $\mathcal{C}$ with at least six vertices.


\begin{enumerate}
\item $G \in \mathcal{C}_=$. If all  segments of $G$ have two vertices,  then the result is trivial since $Ls(G) =  0 $ if $G$ contains an even number of segments,  and $Ls(G) =  2$ if $G$ contains an odd number of segments.

Otherwise, to prove that $Ls(G)    \geq 0 $, note that Left is able to play a $4$-move. Indeed, let $2j$ be the even endpoint of a segment with at least four vertices. Then, playing $2j+2$ (if the even endpoint is the smallest point of the segment) or playing $2j-2$ (otherwise) is a $4$-move. Let $2i$ be this $4$-move.

The resulting relevant  directed graph $G_{2i} = G \setminus \mbox{Rmv}(G_{2i})$ induced by the vertex $2i$ is an element of $ \mathcal{C}_=$.  By the induction hypothesis, we have $Rs(G_{2i}) \in  \{0, -2, -4\}$.  Thus, we obtain
$$Ls(G)  \geq  4 +  Rs(G_{2i})  \geq 0. $$

In other words, Left can obtain a relative score of at least $4-4 = 0$ by playing a 4-move in $G$.

To show that $Ls(G)\leq 4$, we explore all possibilities for the first move of Left. Let $2i$ be the vertex played by Left.
\begin{itemize}
\item If the move at $2i$ is a $4$-move, we know that $Rs(G_{2i}) \in  \{0, -2, -4\}$ from which it follows that the final relative score of Left in $G$  is at most  $4+Rs(G_{2i})  \leq 4 $.

The same argument holds if Left plays a  $2$-move: the final relative score of Left in $G$  is at most  $2+Rs(G_{2i})  \leq 2 $.
\item If the move at $2i$ is a  $3$-move then the resulting relevant  directed graph $G_{2i}$ is such that $G_{2i} \in  \mathcal{C}_+$ since Left has necessarily removed one vertex in $L$ and two vertices in $R$. Thus, by the induction hypothesis, we have $Rs(G_{2i}) \in  \{1, -1, -3 \}$. Thus, if $2i$ is chosen by Left, then  the final relative score of Left in $G$  is $3+Rs(G_{2i})  \leq 4 $.
\end{itemize}
Because there is no other possible move, it follows that $Ls(G) \leq 4$.

We have shown that  $0 \leq Ls(G) \leq 4  $. From Remark \ref{parity}, we obtain $Ls(G) \in \{0, 2, 4\}$, as desired. Moreover, since any element $G$ of $\mathcal{C}_=$ is its own negative, we have $Rs(G)=-Ls(-G)=-Ls(G)$.

\item  $G \in \mathcal{C}_+$.  If all segments of $G$ have at most three vertices,  then the result is trivial.
  If all segments of $G$ have two vertices except one with five vertices, the result is also trivial.
Otherwise, there is a vertex $2i$ associated with Left which induces a $4$-move since there is an endpoint of even label in a component with at least four vertices (but not five).
The resulting  directed graph $G_{2i}$  is an element of $ \mathcal{C}_+$.  By the induction hypothesis, we have $Rs(G_{2i}) \in  \{1, -1, -3\}$.  Therefore, we obtain
$$Ls(G)  \geq  4 +  Rs(G_{2i})  \geq 1. $$

To show that $Ls(G)\leq 5$, we explore all possibilities. Let $2i$ be the vertex played by Left.

\begin{itemize}
\item If the move at $2i$ is a $4$-move, then by induction, we know that $Rs(G_{2i}) \in  \{1, -1, -3\}$. The final relative score in $G$  is at most  $4+Rs(G_{2i})  \leq 5$.

The same argument holds if  Left plays a  $2$-move. The final relative score of Left in $G$  is at most  $2+Rs(G_{2i})  \leq  3 \leq 5 $.
\item If the move at $2i$ is a $3$-move, we have several subcases:
 \begin{itemize}
\item The odd segment is of size 3 and $2i$ lies in this segment. It results that $G_{2i} \in \mathcal{C}_{=}$. By the induction hypothesis, $Rs(G_{2i}) \in \{0, -2, -4\}$ so that the final relative score of Left is at most $ 3 +  Rs(G_{2i}) \leq 3. $
\item The odd segment is of size  at least 5. In this case, $G_{2i}$ contains two segments of  $ \mathcal{C}_+$ and Right can reply by a border 3-move, say on the vertex $2j + 1$, on such a segment. The resulting relevant  directed graph  $(G_{2i})_{2j+1}$ belongs to $ \mathcal{C}_+$. By the induction hypothesis, $Ls((G_{2i})_{2j+1}) \in \{1, 3, 5\}$.
It follows that final relative score of Left in $G$ is $3 - 3 +  Ls((G_{2i})_{2j+1}) \leq 5$.
\end{itemize}

\item If a vertex $2i$ associated with $L$ induces a $5$-move, the resulting relevant directed graph $G_{2i}$ belongs to $ \mathcal{C}_=$. By the induction hypothesis $Rs(G_{2i}) \in\{0, -2, -4\}$. Consequently, the final relative score of Left given that $2i$ is chosen  in $G$  is at most  $5 + R(G_{2i})  \geq  5 $.
 \end{itemize}

From the above case-by-case analysis, we conclude that $Ls(G)\leq 5$. By Remark \ref{parity}, we get that $Ls(G) \in \{1, 3, 5\}$. For $Rs(G)$, we have the same property as before, as $-G$ belongs to $\mathcal{C}_{-}$, which explains why the values of $Ls$ and $Rs$ are opposite for $\mathcal{C}_{-}$ and $\mathcal{C}_{+}$.

\item  $G \in \mathcal{C}_-$.  If all segments of $G$ have at most three vertices,  then the result is straightforward.

Otherwise, there exists a vertex $2i$ that induces a border $3$-move  on the odd segment belonging to $\mathcal{C}_{-}$. The  resulting relevant  directed graph $G_{2i}$  is an element of $ \mathcal{C}_=$.  By the induction hypothesis, we have $Rs(G_{2i}) \in  \{0, -2, -4\}$.  Thus, we get $$Ls(G)  \geq  3 +  Rs(G_{2i})  \geq   3-4 = -1 .$$

To show that $Ls(G)\leq 3$, we explore all possibilities. Let $2i$ be the vertex played by Left in an optimal strategy.

\begin{itemize}
\item If the move at $2i$ is a $3$-move on the odd segment, then the resulting relevant  directed graph $G_{2i}$  is element of $ \mathcal{C}_=$.  By the induction hypothesis, we have $Rs(G_{2i}) \in  \{0, -2, -4\}$.  Thus, the final relative score of Left is at most $ 3 +  Rs(G_{2i})  \leq  3 $.
\item If the move at $2i$ is a $4$-move, then the resulting relevant  directed graph $G_{2i}$  is an element of $ \mathcal{C}_-$. By the induction hypothesis, we have $Rs(G_{2i}) \in  \{-1, -3, -5\}$. Thus, the final relative score of Left is at most
   $4 + Rs(G_{2i})  \leq 4 -1 = 3$.

The same argument holds if Left plays a $2$-move. The final relative score of Left in $G$  is at most  $2 + Rs(G_{2i}) \leq 3 $.

\item If the move at $2i$ is a $3$-move in an even segment, the resulting  directed graph $G_{2i}$ contains a segment in $ \mathcal{C}_+$.
Then, Right can play a border 3-move, say $2j+1$, on this segment. The resulting relevant  directed graph $(G_{2i})_{2j+1}$  is an element of  $\mathcal{C}_-$.  Thus, by the induction  hypothesis,  $Ls((G_{2i})_{2j+1}) \in \{-1, 1, 3\}$. It follows that the final relative score of Left in $G$  is at most  $3-3 + Ls(G_{2i})_{2j+1} \leq 3-3+3 = 3 $.
\end{itemize}
 The above case-by case analysis ensures that $Ls(G)\leq 3$. By Remark \ref{parity}, we get that $Ls(G) \in \{-1, 1, 3\}$.
\end{enumerate}
\end{preuve}

The given sets are the tightest ones.  By an exhaustive  analysis, we have found that the  directed graph $G$,  element of   $\mathcal{C}_-$  formed by the sum of two segments, one  with 10 vertices,
the other one with 17 vertices, satisfies $Ls(G) = -1$. The other values are reached on small examples having a unique segment, and will be illustrated in the following section. Finally, note that the results of Proposition~\ref{prop:segments} do not hold when there are more odd segments. For example, the sum of two segments of $\mathcal{C}_+$ of sizes $5$ and $17$ has a left score of $6$.

\subsection{Values reached for  simple segments}

We have found some minimal sets of values for some families of sums of segments.
We now improve this result by restricting our study to segments, which allows to remove some lowest values of the previous proposition. In particular, we will show that for a segment $G$, we always have $Ls(G)\geq 1$ and $Rs(G)\leq -1$.

\begin{theorem}
\label{prop:simple}
Let $G \in  \mathcal{C}^{segment}$:
 \begin{enumerate}
\item   If $G \in \mathcal{C}_=$,  then  $Ls(G) \in \{2, 4\}$ and $Rs(G) \in \{-2, -4\}$.
\item   If $G \in \mathcal{C}_{-}$, then  $Ls(G) \in \{1, 3\}$  and $Rs(G) \in \{-1, -3, -5\}$. 

Moreover,  if $\vert L\vert$ is odd, then $Ls(G) = 3$.
\item  If If $G \in \mathcal{C}_{+}$,  then   $Ls(G) \in \{1, 3, 5\}$ and $Rs(G) \in \{ -1, -3\}$.

Moreover,  If $\vert L\vert$ is even, then $Rs(G) = -3$.
\end{enumerate}
\end{theorem}

As for Proposition \ref{prop:segments}, one can easily prove that the given sets are tight, as each value is reached.

\begin{corollary}
 For each $G \in  \mathcal{C}^{segment}$, we have  $Ls(G)  > 0$.
\end{corollary}

If we assume that the winner of  the game is the player which has the the largest score,
the corollary above claims that, for any segment, the first player is always   the winner (and there is no tie).

The proof of Theorem \ref{prop:simple} is not a direct corollary of Proposition \ref{prop:segments}.
We need two preliminary successive lemmas.

\begin{lemma}\label{lem:jonction}
Let $G$ be an element of $ \mathcal{C}$ and let $i$ be an integer such that $2i-1, 2i, 2i+1, 2i+2$ are elements of $V(G)$. 
Let $G'$ be the directed subgraph of $G$ such that $V(G') = \{j \in V(G): j \leq 2i+1 \}$   and 
$G''$ be the directed subgraph of $G$ such that $V(G'') = \{j \in V(G):  j \geq 2i+2 \}$. We have: 
$$  Rs(G) \geq  Rs(G') - \vert V(G'') \vert . $$
\end{lemma}

\begin{preuve} First note the following fact. Because $Rs(G)  =2s_2^L(G)- (\vert V(G') \vert + \vert V(G'') \vert)$  and  $Rs(G') =  2s_2^L(G') - \vert V(G') \vert $, the inequality  $  Rs(G) \geq  Rs(G') - \vert V(G'') \vert  $ can be rewritten as  $s_2^L(G) \geq s_2^L(G')$. Thus, it remains to show that $s_2^L(G) \geq s_2^L(G')$.  

Next, assume that Left  is the second player and  uses the strategy consisting in focusing in $G'$ and playing  in an optimal way in $G'$ (while the game $G'$ is not over) without taking into account $G''$. 
More precisely,  when it is Left turn, Left considers the moves previously played as if they were either moves of the game in $G'$ or ``passing moves'', that is for  moves in $G$  which are not feasible in the game starting from $G'$.

For $k \geq 0$, let $H_k$ be the directed graph obtained after $k$ moves where Left plays the above strategy starting from 
$G$;   and let $H'_k$ be the directed graph obtained after the same $k$ moves, starting from $G'$ with the above convention.
 We will compare the sequences of directed subgraphs $(H'_k)_k$ of $G'$ with the sequences of directed subgraphs $(H_k)_k$ of $G$. 
Notice that, at the origin,  $V(H'_0) =  V(H_0) \cap V(G' ) $.  
\begin{enumerate}
\item If $V(H'_k) =   V(H_k) \cap V(G' ) $ and $2i+1 \in V(H'_k)$, then $2i \in V(H'_k)$ since $H'_k$ is a relevant  directed graph.  Thus, if $2i+3  \in H_k$, then  $2i+1 \notin \mbox{Rmv}(H_k, 2i+3)$. It follows that $V(H'_{k+1}) =   V(H_{k+1}) \cap V(G' ) $ and $2i+1 \in V(H'_{k+1})$ except if the following condition holds: 
$2i+2 \in V(H_k) $ and the $(k+1)^{st}$ move is done  by  Right  at vertex $2i-1$  in the game $G'$, which makes $\{2i-1, 2i,  2i+1\} \subseteq  \mbox{Rmv}(H'_k, 2i-1)$.  In this case,   $2i+1 \notin  V(H'_{k+1})$ and   $V(H'_{k+1}) \cup \{2i+1\} =   V(H_{k+1}) \cap V(G' )$. 
This case is considered in the following point. 

\item  { If} vertices $2i-1, 2i, 2i+1$ are not elements of $V(H'_{k})$ and  $V(H'_{k}) \cup \{2i+1\} =   V(H_{k}) \cap V(G' )$, then, since $H_k$ is a relevant  directed graph, it holds that $2i+2 \in V(H_k)$. Thus, $2i-1, 2i, 2i+1$ do not belong to $V(H'_{k+1})$ and   $V(H'_{k+1}) \cup \{2i+1\} =   V(H_{k+1}) \cap V(G' )$, except if one of the following conditions holds:
\begin{itemize}
\item the $(k+1)^{st}$ move is done by  Right    at $2i+1$;
\item the vertex $2i+3 \in V(H_k) $ and the $(k+1)^{st}$ move is done by  Right   at  $2i+3$. 
\end{itemize}
If one of these alternatives holds,  $V(H'_{k+1}) =   V(H_{k+1}) \cap V(G' )$ and, of course, $2i+1 \notin V(H'_{k+1}) $. 
This case  is considered in the following  point. 
\item { If} $V(H'_{k}) =   V(H_{k}) \cap V(G' )$ and $2i+1 \notin V(H'_{k}) $, then  obviously $V(H'_{k+1}) =   V(H_{k+1}) \cap V(G' )$ and $2i+1 \notin V(H'_{k+1}) $. 
\end{enumerate}

It follows  that,  for each $k$ such that $V(H'_{k }) \neq  \emptyset$, one of the cases described in (1), (2) or (3) holds. Thus $V(H'_{k}) \cap L=   V(H_{k}) \cap V(G' ) \cap L$, which ensures that the strategy described for the Left payer is well defined. Furthermore, for each $x \in H'_{k} \cap L$,   $\mbox{Rmv}(H_k, x) \cap V(G' ) =  \mbox{Rmv}(H'_k, x) $. Let $(x_1, x_2, \ldots) $ be the sequence of moves generated by  the Left player in  $G'$, meaning that the move $x_\ell$ is initiated  in $H'_{2 \ell-1}$. We get:
$$  s_2^L(G') =   \sum_\ell   \vert   \mbox{Rmv}(H'_{2 \ell-1}, \ell)  \vert  =   \sum_\ell   \vert   \mbox{Rmv}(H_{2 \ell-1}, x_\ell) \cap V(G' ) \vert  \leq  \sum_\ell   \vert   \mbox{Rmv}(H_{2 \ell-1}, x_\ell) \vert =  s_2^L(G), $$
which is the expected result. 
\end{preuve}

\begin{lemma}\label{lem:twoseg}
Let  $G \in \mathcal{C}^{}$  such that $G  =S_{s_1} + S_{s_2}$, where   for $i \in \{1, 2\}$,  $S_{s_i}$ is a segment with $s_i$  vertices. \begin{enumerate}
\item If $s_1 = s_2$  and   $G \in \mathcal{C}_=$, then  $Rs(G) = 0$.
\item If $s_1 =  s_2 $,  $S_{s_1} \in  \mathcal{C}_-$ and  $S_{s_2} \in  \mathcal{C}_+ $, then $Rs(G) = 0$.
\item If $s_1  = 2k$ and  $s_2  = 2k+2$,   then  $Rs(G) = -2$.
\item If $s_1  = 2k$,  $s_2  = 2k+1$, and $G \in  \mathcal{C}_+ $, then  $Rs(G) = -1$.
\item If $s_1  = 2k$,  $s_2  = 2k-1$, and  $G  \in  \mathcal{C}_+ $,  then  $Rs(G) \geq -1$.
 \end{enumerate}
\end{lemma}

\begin{preuve} Pick any $G = G  =S_{s_1} + S_{s_2}$ as hypothesized.
\begin{enumerate}
\item If  $H$ is a segment of $\mathcal{C}_= $, then $H = -H$.  This allows to apply  Corollary \ref{cor:G-G} where $H= S_{s_1}$ in such a way that $G  = H  + (-H)$. 
\item If  $H$ is a segment of $\mathcal{C}_+ $, then $-H$ is the segment of $\mathcal{C}_- $ with the same number of vertices. This allows to apply  Corollary \ref{cor:G-G},  where $H= S_{s_1}$ in such a way that $G  = H  + (-H)$. 
\item Here, we proceed by induction on $k$. For $k = 1$, the  result is obvious.
Now assume that  the result is true for a $k \geq 1$ and take  $G = S_{s_1}+S_{s_2}$, such that $S_{s_1}$ is a segment with $2(k +1)$ vertices and $S_{s_2}$ is a segment with $2(k +1) +2$ vertices.
Let $y$ denote the $R$-endpoint of $S_{s_2}$, $x$  be the unique vertex such that $(x, y)$ is a directed arc of $S_{s_2}$  and $S_{2} $ be the instance containing only the vertices $x$ and $y$ and the arc $(x, y)$. Let $G' = G \setminus \{x, y\}$.  Hence, $G'$ is the sum  of two segments with $2k +2$ vertices for each segment. Thus, by point (1),  $Rs(G') = Ls(G')  = 0$.

Furthermore,  $S_{2}$ and $G'$  allow  to apply Lemma   \ref{lem:jonction} with $G'' = S_2$, so that  we get
$$  Rs(G) \geq  Rs(G') - \vert V(S_{2}) \vert  = 0-2   = -2. $$
To prove the converse inequality, assume that the first move of Right is a 4-move at a vertex $y' $ of $S_{s_2}$
and,  then,  players both play  optimally. The relative score of Left with this strategy is $-4 + Ls(G_{y'})$.
The  lengths of the segments constituting $G_{y'}$  are  $2k$ and $2(k+1)$. Thus,  since these lengths are even,  by symmetry, we have  $Ls(G_{y'}) = - Rs(G_{y'})$. Furthermore,   the  induction hypothesis applies to $G_{y'}$ so that $Rs(G_{y'}) = - 2$. We obtain $Rs(G)  \leq -2$. 

\item 
 We use the same procedure as in point (3).
Let $x$ denote an endpoint of $S_{s_2}$. Notice that $x$ is necessarily a $L$-vertex.  We denote by  $S_{1}$  the directed  graph reduced to the  singleton $\{x\}$,  and $G' = G \setminus \{x\}$. We can apply point (1) to  $G'$,  thus  $Rs(G') = Ls(G')  = 0.$ Moreover,  Lemma   \ref{lem:jonction} applies  with $G'' = S_1$ (except when $k = 1$, which is a trivial  case), in order to obtain $$  Rs(G) \geq  Rs(G') - \vert V(S_{1}) \vert  = 0-1  = -1. $$
To prove the converse inequality, note that Right can obtain a relative score equal to $-1$ by playing a border 3-move $y'$ on $S_{s_2}$. Indeed, by point (3) applied to the  directed graph $G_{y'}$, we get  $Ls(G_{y'}) = 2$ and thus a relative score of $-1$ for Right.

\item  The proof is similar to point (4). Let $x$ denote the endpoint of $S_{s_1}$ which is  a $L$-vertex.  We denote by $S_{1}$  the  directed graph reduced to the singleton $\{x \}$ and $G' = G \setminus \{x\}$. We can apply point (2) to $G'$,  and so  $Rs(G') = Ls(G')  = 0.$
Furthermore,  Lemma   \ref{lem:jonction} applies with $G'' = S_{1}$, thus we get $$  Rs(G) \geq  Rs(G') - \vert V(S_{s_1}) \vert  = 0-1  = -1. $$
\end{enumerate}
\end{preuve}

We now have the material  to prove Theorem \ref{prop:simple}.

\begin{preuve}[of Theorem \ref{prop:simple}].
The proof  is done by ``cutting the segment in quasi-equal parts''  as follows.

\begin{enumerate}
\item Except in the trivial cases of size  at most 6,  Left can play a cutting move $x = 2i$  in such a way that the resulting  directed graph $G_{x}$ is formed by two segments satisfying the hypothesis of point 4 of Lemma \ref{lem:twoseg}. Thus, after playing the cutting move, the final relative score of Left is $3 +Rs(G_{x})  \geq  3-1 = 2$. This ensures that $Ls(G) \geq 2$. From Proposition \ref{prop:segments}, we have  $Ls(G) \in \{0, 2, 4\} $, which gives that  $Ls(G) \in \{2, 4\} $. By symmetry,  we get the result for  $Rs(G)$.

 \item Except in the trivial cases of size  at most 5,  Left can play a cutting move $x$ in such a way that the resulting  directed graph $G_{x}$ is formed by two segments, each of them containing an even number of vertices, satisfying the hypothesis of either point 1 or point 3.  of Lemma \ref{lem:twoseg} according to the value  $\vert L\vert +  \vert R \vert$ modulo $4$. Precisely,
 \begin{itemize}
\item when $\vert L\vert +  \vert R \vert  \equiv 3 \mod [4]$,  then point 1 of Lemma \ref{lem:twoseg} applies.  Thus, after playing the cutting move,
the final relative
 score of Left is $3 +Rs(G_{x})  =   3+0 = 3$. This ensures that $Ls(G) \geq 3$. From Proposition \ref{prop:segments}, we have  $Ls(G) \in \{-1, 1, 3\} $, which gives that   $Ls(G) =  3$;
\item  when $\vert L\vert +  \vert R \vert  \equiv 1 \mod [4]$,  then point 3  of Lemma \ref{lem:twoseg} applies.   Thus, after playing the cutting move,
the final relative
 score of Left is $3 +Rs(G_{x})  \geq   3-2 = 1$. This ensures that $Ls(G) \geq 1$. From Proposition \ref{prop:segments}, we have  $Ls(G) \in \{-1, 1, 3\} $ from which il follows that   $Ls(G) =  \{1, 3\}$.
 \end{itemize}
\item This case is symmetric to the previous  case.
\end{enumerate}
  \end{preuve}

Concerning the case of cycles, the result follows from the case of segments since all the first moves are equivalent moves and lead to a segment.

\begin{corollary}
 \label{cor:cycle}
Let $G \in  \mathcal{C}^{cycle}$.
We have  $Ls(G) \in \{0, 2\}$.
Moreover,
 if $\vert L\vert$ is odd, then $Ls(G) = 0$.
\end{corollary}

 \begin{preuve}
Let $i$ be the vertex where the first move is initiated by Left.
Note that  $Ls(G) = 3 - Rs(G_i)$,  and  $G_i$ satisfies point 3 of   Theorem \ref{prop:simple} from which the result follows.
 \end{preuve}

\subsection{The score sequence of segments }

For a segment of size $n$, the next question that arises is about the complexity of the computation of the exact scores $Ls$ and $Rs$. More precisely, since from Theorem~\ref{prop:simple} there is a limited set of values for them, one can wonder whether there is a kind of regularity for the sequences $(Ls(S_n))_{n\geq 1}$ and $(Rs(S_n))_{n\geq 1}$, where $S_n$ is the segment of size $n$ with the convention that, when $n$ is odd, we consider the segment $S_n$ to be in $\mathcal{C}_-$.  As it is the case for some combinatorial games played on strings (such as subtraction games or octal games, see~\cite{Siegel}), a regular behavior of the sequence is a key to get a polynomial time algorithm that computes the value (or score) of the game. \\

Table~\ref{tab:segments} summarizes the first values of these two sequences (computed by our program).

\begin{table}[htbp]
{\small
$$\begin{array}{|c||r|r|r|r|r|r|r|r|r|r|r|r|r|r|r|r|r|r|r|}
\hline
n & 1& 2& 3& 4& 5& 6& 7& 8& 9& 10& 11& 12& 13& 14& 15 & 16 & 17 & 18 & 19\\
\hline
Ls(S_n)&-1&2 & 3 & 4 & 1 & 2 & 3 & 2 & 1 & 2 & 3 & 2 & 1 & 4 & 3 & 2 &1&2&3\\
Rs(S_n)&-1&-2 & -3 & -4 & -5 & -2 & -1 & -2 & -3 & -2 & -1 & -2 & -3 & -4 & -3 & -2& -3&-2&-3\\
\hline
\hline
n & 20& 21& 22& 23& 24& 25&26& 27& 28& 29& 30& 31& 32& 33& 34& 35 & 36 & 37 & 38 \\
\hline
Ls(S_n)&2&1 &  4 & 3 & 2 & 1 & 2 & 3 & 2 & 1 & 4 & 3 & 2 & 1 & 2 &3&2&1&2\\
Rs(S_n)&-2&-5  & -4 & -3 & -2 & -3 & -2 & -3 & -2 & -5 & -4 & -3 & -2 & -3 & -2& -3&-2&-5&-2\\
\hline
\end{array}
$$
}
    \caption{First values of the sequences $(Ls(S_n))_{n\ge 1}$ and $(Rs(S_n))_{n\ge 1}$.}
    \label{tab:segments}
\end{table}

Some regularities can be found in the sequences. In particular, Theorem \ref{prop:simple} proves that $Ls(S_n)=3$ for $n\equiv 1 \bmod[4]$ and $n\neq 5$. But the sequences have been computed up to $n=80$ and no periodicity appears yet. In particular, between the sizes $38$ and $76$, both sequences are periodic with period of size $4$. Unfortunately, we have $Rs(S_{77})=-5$, that breaks the previous regular behavior. However, as for octal games, a periodicity conjecture can be established in view of the very restricted number of possible scores. 

\begin{conjecture}
The sequences $Ls(S_n)$ and $Rs(S_n)$ are ultimately periodic.
\end{conjecture} 

If a proof of this conjecture is currently out of reach because of the global chaotic behavior of the sequences, the sub-regularities that can be found lead to the following more affordable conjectures.

\begin{conjecture}
For a segment $S_n$ of $\mathcal{C}_-$ with $n\equiv 1\bmod 4$, we have $Ls(S_n)=1$.
\end{conjecture} 

\begin{conjecture}\label{conj:2mod4}
  For a segment $S_n$ with $n\equiv 0\bmod 4$ and $n>4$, we have $Ls(S_n)=2$.
\end{conjecture}

Actually, the lack of periodicity is due to the ``rare'' scores $(Ls(S_n),Rs(S_n))\in\{(4,-4),(1,-5)\}$ for some values of $n$ that we did not manage to characterize. In particular, the last score $(4,-4)$ has been found for $n=30$, and one can wonder whether it is the last one. In that case, Conjecture~\ref{conj:2mod4} could be extended to all the even values of $n$ greater than $30$.

\section{Position versus  Number}

A natural question arises about the game  \Jeu: Is it better  for a player to have initially a larger number of vertices, or to have few vertices, but in a strategical position? The intuition is that some strategical sites are very important. We quantify  this intuition  (see Theorem  \ref{th:smallbig}) by giving  instances $G$ where the proportion of $L$-vertices is as small as desired  and the proportion of vertices won by Left at the end of the game is as closed to 1 as desired. We then prove that this is not true for oriented paths (see Theorem \ref{th:proppath}).

\subsection{An instance with a large number of $R$-vertices but good for Left}
  
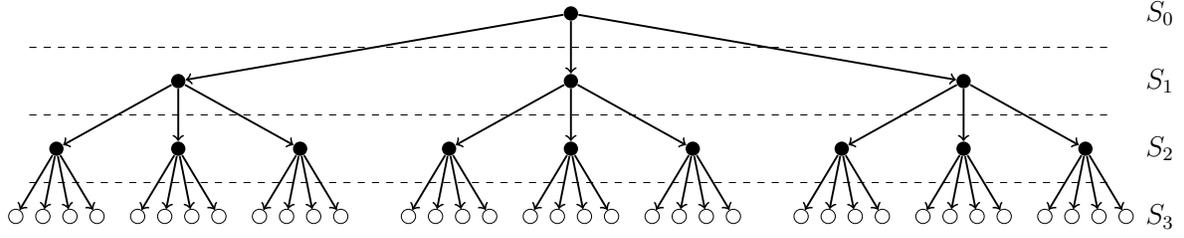
\begin{figure}[h]
  \begin{center}
    \scalebox{0.9}{\begin{tikzpicture}
      \node[vertex] (0) at (0,0) {};

     \foreach \K in {-1,0,1}{
      \begin{scope}[shift={(\K*5.8,-1)}]
     \node[vertex] (1\K) at (0,0) {};
     \draw[edge] (0) to[->] (1\K);
     \foreach \J in {-1,0, 1}{
  \begin{scope}[shift={(\J*1.8,-1)}]
       \node[vertex] (2\K\J) at (0,0) {};
           \draw[edge] (1\K) to[->] (2\K\J);
 
           \foreach \I in {0,1,2,3}{
             \node[vertexW] (3\K\J\I) at (-0.6+\I*0.4,-1) {};
             \draw[edge] (2\K\J) to[->] (3\K\J\I);}
      
     \end{scope}}
     \end{scope}   }

     \draw[dashed] (-8,-0.5) -- (8,-0.5) ; 
     \draw[dashed] (-8,-1.5) -- (8,-1.5) ; 
     \draw[dashed] (-8,-2.5) -- (8,-2.5)  ;
     \draw (8.7,0) node {$S_0$}; 
     \draw (8.7,-1) node {$S_1$}; 
     \draw (8.7,-2) node {$S_2$}; 
     \draw (8.7,-3) node {$S_3$}; 
    
    \end{tikzpicture}}
  \end{center}
  \label{fig:Tnc}
  \caption{The oriented tree $T_2^4$.}
  \end{figure}

We define our instances as follows. For each pair $n$, $c$ of integers, with $n \geq 0$ and  $c \geq 1$, let $T_n^c$ be the oriented tree of depth $n+1$,
such  that all vertices of level $k$, $0 \leq k \leq  n-1$ have exactly 3 successors and
vertices of level $n$  have exactly $c$ successors. 
We denote by $S_i$ the vertices of level $i$. 
The set of $R$-vertices is the set of leaves at level $S_{n+1}$, whereas the other vertices are in $L$, $L=\cup_{k= 0}^n S_k$.
See Figure \ref{fig:Tnc} for a representation of $T_2^4$.

For $0 \leq k \leq  n$,  we have  $ \vert S_k \vert  = 3^k$   and $ \vert S_{n+1}  \vert  = c \, 3^n$. Thus
$$ \vert V(T_n^c) \vert  = 1 + 3 + \cdots + 3^n + c\, 3^n  =  \frac{3^{n+1}-1 }{2 } + c \ 3^n. $$

\begin{theorem} \label{th:smallbig}
Let   $J^c_n$  be the instance    $T_n^c   +  T_n^c $, formed by two distinct copies of  $T_n^c $.
For each  $\epsilon >0$,  for $n$ and $c$ sufficiently large,  we have:
 $$ \frac{\vert L(J^c_n) \vert  }{ \vert V(J^c_n)\vert } \leq \epsilon \quad   \mbox{ and } \quad  \frac{s_2^L(J^c_n) }{\vert V(J^c_n) \vert } \geq 1-  \epsilon. $$
\end{theorem}

Note that,  from Corollary \ref{cor:nopass}, if $$ \frac{\vert  s_2^L(J^c_n)\vert }{\vert  V(J^c_n) \vert } \geq 1-  \epsilon,$$ 
then  we also have $$\frac{\vert s_1^L(J^c_n)\vert }{\vert V(J^c_n)\vert } \geq 1-  \epsilon.$$\\

For  the proof  of  the theorem  above, we first need a lower bound on $s_1^R(J^c_{n})$.

\begin{lemma} \label{lem:T_n}
Let $y$ be  an  $R$-vertex  of an instance  $T_n^c$. Then  $(T_n^c)_y=J^c_0  +  J^c_1 +  \ldots + J^c_{n-1} $.
\end{lemma}

\begin{preuve}
For $0 \leq k \leq n$,  $\mbox{Rmv}(T_n^c, y) \cap S_k$  is a singleton, say $x_k$,
and $\mbox{Rmv}(T_n^c, y) \cap S_{n+1} = \{z :(p(y), z) \in A^c_n\}$ where $p(y)$ is the predecessor of $y$.
Let $S$ be the set of vertices of $(T_n^c)_y$  with no predecessor.  For  $1 \leq k \leq n$, the set
$S \cap S_k$ is a pair,  say $\{z_k, z'_k \}$.
The  directed graph induced by the successors of $z_k$ (respectively $z'_k$) is identical to the   directed graph  $T_{n-k}^c$, from which we get the result.
\end {preuve}

Using the above lemma, we obtain a {recurrence} relation for $s_1^R(J^c_n)$:

\begin{lemma}\label{lem:recurrence}
  Let $n\geq 0$, we have:
$$ s_1^R(J^c_{n}) \leq n+1+c +  \sum_{i=0}^{n-1} s_1^R(J^c_i).$$
\end{lemma}

\begin{preuve}
 The first move of Right in $J^c_n$ is necessarily a leaf $y$ of $T^c_n$, leading, by Lemma \ref{lem:T_n} to the  directed graph $G=J^c_0+J^C_1+...+J^c_{n-1}+T_n^c$.
  Thus,
  $$ s_1^R(J^c_{n})= n+1+c +  s_2^R(G).$$

  By definition of the score, $s^R_2(G)=\min_{x \in L} \bigl\{s^R_1(G_x) \bigr\}$. Let $x\in L(G)$ be the root of $T_n^c$. If Left plays $x$, then she takes all the component $T_n^c$. We thus obtain: 
  $$s^R_2(G)\leq  s^R_1(G_x) =  s^R_1(J^c_0+J^C_1+...+J^c_{n-1}). $$
   By Corollary \ref{cor:milnor}, $s^R_1(J^c_0+J^C_1+...+J^c_{n-1})\leq \sum_{i=0}^{n-1} s^R_1(J^c_i)$,  which concludes the proof.
   \end{preuve}

Applying inductively Lemma \ref{lem:recurrence}, we can upperbound $s_1^R(J^c_n)$:

\begin{corollary}\label{cor:ubsR}
    Let $n\geq 0$, we have: $s_1^R(J^c_{n}) \leq 2^n(n+c)+1.$
\end{corollary}

\begin{preuve}
  We prove by induction on $k$ that the following inequality holds for any $0\leq k\leq n$:
\begin{equation}\label{eq:rec}
  s_1^R(J^c_n)\leq n+1+c+\sum_{i=0}^{k-1}(n+c-i)2^i+ 2^{k}\sum_{i=0}^{n-1-k} s_1^R(J^c_i).
\end{equation}

The initialization $k=0$ is exactly Lemma \ref{lem:recurrence}.
Assume the result is true for $k$ such that $0\leq k <n$.
We apply Lemma \ref{lem:recurrence} to $s_1^R(J^c_{n-1-k})$ and obtain:
\begin{equation*}
  \begin{split}
    s_1^R(J^c_n)&\leq n+1+c+\sum_{i=0}^{k-1}(n+c-i)2^i+ 2^k\left( n+c-k+\sum_{i=0}^{n-k-2}s_1^R(J^c_i)\right)+ 2^{k}\sum_{i=0}^{n-2-k} s_1^R(J^c_i)\\
    &\leq  n+1+c+\sum_{i=0}^{k}(n+c-i)2^i+ 2^{k+1}\sum_{i=0}^{n-2-k} s_1^R(J^c_i)
  \end{split}
  \end{equation*}

which corresponds to the relation (\ref{eq:rec}) for $k+1$.

Equation (\ref{eq:rec}) with $k=n$ gives us:
\begin{equation*}
  \begin{split}
    s_1^R(J^c_n)&\leq n+1+c+\sum_{i=0}^{n-1}(n+c-i)2^i\\
    &\leq  n+1+c+(n+c)\sum_{i=0}^{n-1}2^i \\
    &\leq 2^n(n+c)+1.
  \end{split}
  \end{equation*}
\end{preuve}

We are now able to prove Theorem \ref{th:smallbig}.

\begin{preuve}  [(of Theorem \ref{th:smallbig})]
  First, we have $\vert L(J^c_n) \vert = 2 (1 + 3 + \ldots + 3^n) = 3^{n+1} - 1$ and   
$$\vert V(J^c_n) \vert =  3^{n+1} - 1 + 2c 3^n \geq   (3^{n+1} - 1 )\left( 1 +   \frac{2c} {3}\right)$$

Thus, for   $$\frac{1}{1 +   \displaystyle\frac{2c} {3}} \leq {\epsilon},\quad  \mbox{or equivalently} \quad c  \geq \frac{3}{2} \left(\frac{1}{\epsilon } -1\right),$$ we have $$\frac{\vert L(J^c_n) \vert  }{ \vert V(J^c_n)\vert } \leq \epsilon.$$\\

On the other hand, we have, by definition of the score, $s_2^L(J^c_n)=|V(J_n^c)|-s_1^R(J^c_n)$.
Using the upper bound of Corollary \ref{cor:ubsR}, we get:

\begin{equation*}
  \begin{split}
    \frac{s_2^L(J^c_n)}{|V(J_n^c)|}&= 1-\frac{s_1^R(J^c_n)}{|V(J_n^c)|}\\
    &\geq 1 - \frac{2^n(n+c)+1}{3^{n+1}-1+2c3^n}
  \end{split}
\end{equation*}

Finally,
$$\lim_{n\to +\infty} \frac{2^n(n+c)+1}{3^{n+1}-1+2c3^n}=0 \quad \text{and thus} \quad \lim_{n\to +\infty} \frac{s_2^L(J^c_n)}{|V(J_n^c)|}=1.$$

Hence,  whenever the integer $c$ has been chosen in order to have  $ \frac{\vert L(J^c_n) \vert  }{ \vert V(J^c_n)\vert } \leq \epsilon $, then the integer $ n$ can be chosen in order to have $\frac {s_2^L(J^c_{n})}{ \vert J^c_{n} \vert } \geq 1 - \epsilon$.
\end{preuve}

\subsection{Guarantee of the score in quasi-paths}

In the previous example, some $L$-vertices have a large number of successors whereas the $R$-vertices have at most $n+1$ predecessors, making possible for Left to win a lot of vertices in a single move. We now introduce  a class of  directed graphs where is not true anymore. It follows that   Left  cannot win too many vertices compared to the initial situation.

\begin{definition}
Let $G = (L \cup R, A)$ be relevant graph of $\mathcal{S}$. The directed graph $G$ is a quasi-path if it is an alternative sequence of vertices of $L \cup R$ and arcs of $A$ that forms a path in the underlying undirected graph of $G$. 
Let $\mathcal{D} \subseteq \mathcal{S}$ be the subclass of relevant directed graphs such that each $G \in \mathcal{D}$ is a collection of quasi-paths.
\end{definition}

Let $G = (L \cup R, A)$ be a collection of quasi-paths of $\mathcal{D}$. Denote by $n$ the cardinality of $L \cup R$, by  $n_L$ the cardinality of the set of $L$-vertices and by $n_R$ the cardinality of the set of $R$-vertices of $G$. In particular, $u$ stands for the ratio $u=n_L/n$.


\begin{theorem}\label{th:proppath}
  Let $G = (L \cup R, A)$ be a collection of quasi-paths of $\mathcal{D}$. 
 Then, it holds that 
  $$\frac{s_2^L(G)}{|V(G)|}\leq \frac{2+u}{3},$$
  where, as usual, $V(G) = L \cup R$. 
\end{theorem}

Before proceeding to the proof of the statement of Theorem \ref{th:proppath}, consider the following example. Assume that graph $G \in \mathcal{D} $ contains a majority of $R$-vertices. Then, Theorem \ref{th:proppath} indicates that the  ratio of vertices won by Left is at most $$\frac{2+ \displaystyle \frac{1}{2}}{3}  =\frac{5}{6}$$ 
meaning that this ratio cannot be very close to 1. 

\begin{preuve}
Let $G = (L \cup R, A)$ be an element of $\mathcal{D}$, and assume that Right starts the game.
It is enough to prove that there is a strategy for Right where the number of $R$-vertices won by $R$
is at least equal to half the number of $R$-vertices won by Left. Indeed, assume that such a strategy exists. Then, Left can win at most all the $L$-vertices and $2/3$ of the $R$-vertices. Therefore, if such a strategy for $R$ exists,then we obtain the desired result:

  \begin{equation*}
    \begin{split}
      \frac{s_2^L(G)}{n}&\leq \frac{n_L}{n}+\frac{2}{3}\frac{n_R}{n}\\
      &\leq u+\frac{2}{3}(1-u) \\
      &=\frac{2+u}{3}.
    \end{split}
  \end{equation*}

  It remains to prove that this strategy for $R$ exists. 
    Let $x\in L$ be a move induced by Left that maximises the quantity $|\mbox{Rmv}(G,x)\cap R)|$, i.e.  the number of $R$-vertices that Left can obtain in one move. Denote by $k$ this maximal quantity. 
  
Because $G \in \mathcal{D}$, we have  $\mbox{Rmv}(G,x) = \{v_1, v_2, \ldots,v_p \}$, where,  for $1 \leq i < p$, $(v_i, v_{i+1}) \in A$. Let  $j$ (respectively $j'$) be the lowest (respectively the largest) integer such that  $v_j$ (respectively $v_{j'}$) is an element  of  $\mbox{Rmv}(G,x) \cap R$. 
Because each vertex $v_\ell$ of $\mbox{Rmv}(G,x) \cap R$ is  a successor of $x$ in $G$,  $v_\ell$ is necessary localized in a directed path of $G$ from $x$ to either $v_j$ of $v_{j'}$. It follows that 
$$\mbox{Rmv}(G,x)\cap R(G) \subseteq \mbox{Pred}(G,v_j) \cup \mbox{Pred}(G,v_{j'})  \subseteq \mbox{Rmv}(G,v_j) \cup \mbox{Rmv}(G,v_{j'}). $$
 Without loss of generality, assume that $ \vert \mbox{Rmv}(G,v_j)\cap R(G)\vert  \geq \vert \mbox{Rmv}(G,v_j')\cap R(G)\vert $. It follows that  
 $$|\mbox{Rmv}(G,v_j)\cap R(G)|\geq k/2.$$
 Assume now that the first move of Right is to play $v_j$. Then Right wins at least $k/2$ $R$-vertices. Let $x_1$ be the response by Left. 
Because $x_1\notin \mbox{Rmv}(G,y_1)$, by Lemma~\ref{lem:commute}, we have
  $\mbox{Rmv}(G_{y_1},x_1)=\mbox{Rmv}(G,x_1)$, and so, by maximality of $x$, $$|\mbox{Rmv}(G_{y_1},x_1)\cap R|  = |\mbox{Rmv}  (G ,x_1)\cap R|    \leq k. $$
  Thus, in the first two moves, Right wins at least half of the number of $R$-vertices that Left took. The remaining graph $(G_{y_1})_{x_1}$ is still an element of $\mathcal{D}$,  and Right can continue to play his strategy to obtain at the end of the game at least half of the number of $R$-vertices won by Left, concluding the proof.
  \end{preuve}

\end{document}